\documentclass[12pt,centertags]{amsart}
\usepackage{amsmath,amstext,amsthm,a4,amssymb,amscd}
\usepackage[mathscr]{eucal}
\usepackage{mathrsfs}
\usepackage{epsf}
\numberwithin{equation}{section}

\newtheorem{thm}{Theorem}[section]

\newtheorem{prop}[thm]{Proposition}

\theoremstyle{definition}
\newtheorem{rem}[thm]{Remark}
\theoremstyle{definition}

\newcommand{\be}{\begin{eqnarray}}
\newcommand{\ee}{\end{eqnarray}}

\newcommand{\comment}[1]{}

\begin{document}

\title{A note on the Lichnerowicz vanishing theorem for proper  actions}
%\date{May 7, 2009}
%\date{\today}
\author{Weiping Zhang}

\address{Chern Institute of Mathematics \& LPMC, Nankai
University, Tianjin 300071, P.R. China}
\email{weiping@nankai.edu.cn}

\begin{abstract}     We prove a Lichnerowicz type vanishing theorem for non-compact spin manifolds admiting proper cocompact actions. This extends a previous result of Ziran Liu  who proves it for the case where the acting  group is unimodular.
\end{abstract}

\maketitle
%\tableofcontents

\setcounter{section}{-1}
%%%%%%%%%%%%%%%%%%%%%%%%%%%%%%%%%%%%%%%%%%%
%%%%%%%%%%%%%%%%%%%%%%%%%%%%%%%%%%%%%%%%%%%

\section{Introduction} \label{s0}
%%%%%%%%%%%%%%%%%%%%%%%%%%%%%%%%%%%%%%%%%%%

A classical  theorem of Lichnerowicz \cite{L63} states that if an
 even dimensional closed smooth spin manifold    admits a Riemannian metric
of positive scalar curvature, then the index of the associated Dirac operator  vanishes.   
In this note   we prove an extension of this vanishing theorem to the case where a  (possibly non-compact) spin manifold $M$ admiting a proper cocompact  action by a locally compact group $G$. 

To be more precise, recall that for such an action, a so called $G$-invariant index %, denote  ${\rm ind}_GD_+$ 
has been defined by Mathai-Zhang in \cite{MZ08}. Thus it is natural to ask whether this index vanishes if 
$M$ carries a $G$-invariant Riemannian metric of positive scalar curvature.  Such a result has indeed been proved by Liu in  \cite{L13} for the case of unimodular $G$. In this note we extend Liu's result to the   case of general $G$.

We will recall the defnition of the Mathai-Zhang index \cite{MZ08} and state the    main result as Theorem \ref{t2}  in Section \ref{s1};  and then  prove   Theorem \ref{t2}   in Section \ref{s2}.

 %%%%%%%%%%%%%%%%%%%%%%%%%%%%%%%%%%%%%%%%%%%
%%%%%%%%%%%%%%%%%%%%%%%%%%%%%%%%%%%%%%%%%%%
\section{A vanishing theorem for the Mathai-Zhang index}\label{s1}
%%%%%%%%%%%%%%%%%%%%%%%%%%%%%%%%%%%%%%%%%%%
 
   Let $M$ be an    even dimensional spin manifold. Let $G$ be a locally compact group  which acts on $M$ properly and cocompactly, where by proper we mean that the  map
$$G \times M \rightarrow M \times M,\ \ \  (g, x)\mapsto (x, gx),$$
is proper (the pre-image of a compact subset is compact), while by cocompact we mean that the quotient  $M/G$ is compact.
We also assume that $G$ preserves the spin structure on $M$. 

Given a $G$-invariant Riemannian metric $g^{TM}$ (cf.  \cite[(2.3)]{MZ08}), we can construct canonically a $G$-equivariant  Dirac operator $D:\Gamma(S(TM))\rightarrow \Gamma(S(TM))$ (cf. \cite{LM89} and \cite{MZ08}), acting on the Hermitian spinor bundle $S(TM)=S_+(TM)\oplus S_-(TM)$.  Let $D_\pm : \Gamma(S_\pm(TM))\rightarrow \Gamma(S_\mp(TM))$ be the obvious restrictions. 

 Let $\|\cdot \|_0$ be the standard   $L^2$-norm on $\Gamma(S(TM))$, let $\|\cdot \|_1$ be a (fixed) $G$-invariant Sobolev $1$-norm. Let ${\bf H}^0(M,S(TM) )$ be the completion of $\Gamma(S(TM)  )$ under$\|\cdot \|_0$.
Let $\Gamma(S(TM) )^G$ denote the space of $G$-invariant smooth sections of $S(TM) $.

Recall that by the compactness of $M/G$, there exists  a compact subset $Y$ of $M$ such that $G(Y ) = M$ (cf. \cite[Lemma 2.3]{P89}). Let $U,\ U'$ be two open subsets of $M$ such that $Y \subset U$ and that the closures $\overline{U}$ and $\overline{U'}$ are both compact in $M$, and that $\overline U \subset U'$. 
Following \cite{MZ08}, let $f \in  C^\infty(M)$ be a nonnegative function such that $f|_U = 1$ and ${\rm Supp}(f) \subset U'$.

Let ${\bf H}^0_f(M,S(TM) )^G$ and  ${\bf H}^1_f(M,S(TM) )^G$ be the completions of
$\{fs : s \in \Gamma(S(TM)  )^G\}$ under $\|\cdot \|_0$ and $\|\cdot \|_1$ respectively.
Let $P_f$ denote the orthogonal projection from ${\bf H}^0(M,S(TM) )$ to  ${\bf H}^0_f(M,S(TM) )^G$.
Clearly, $ P_fD$ maps ${\bf H}^1_f(M,S(TM) )^G$ into ${\bf H}^0_f(M,S(TM) )^G$. 

We   recall a basic result from \cite[Proposition 2.1]{MZ08}.

\begin{prop} \label{t1}   The operator
$
P_fD : {\bf H}^1_f(M,S(TM) )^G \rightarrow {\bf H}^0_f(M,S(TM) )^G
$
is a Fredholm operator.
\end{prop}

 It has been shown in \cite{MZ08} that ${\rm ind}(P_fD_+)$ is independent of the choice of the cut-off function $f$, as well as the $G$-invariant metric  involved. Following \cite[Definition 2.4]{MZ08}, we denote ${\rm ind} (P_f D_+ )$ by ${\rm ind}_G ( D _+ )$.

The main result of this note can be stated as follows.

\begin{thm}\label{t2} If there is a $G$-invariant metric $g^{TM}$ on $TM$ such that its scalar curvature $k^{TM}$ is positive over $M$, then ${\rm ind}_G ( D _+ )=0$.
\end{thm}

\begin{rem} \label{t3} If $G$ is unimodular, then Theorem \ref{t2} has been proved in \cite{L13}. Our proof of Theorem \ref{t2} combines the method in \cite{L13} with a simple observation that in order to prove the vanishing of the index, one need not restrict to   self-adjoint operators. 
\end{rem}

%%%%%%%%%%%%%%%%%%%%%%%%%%%%%%%%%%%%%%%%%%%
\section{Proof of Theorem \ref{t2}
}\label{s2} 

Following \cite[(2.16)]{MZ08}, let $\widetilde D_{f,\pm}:  {\bf H}^1_f(M,S_\pm(TM) )^G \rightarrow {\bf H}^0_f(M,S_\mp(TM) )^G
$ be defined by that for any $s\in \Gamma(S_\pm(TM) )^G$,
\begin{align}\label{2}
\widetilde D_{f,\pm} (fs)=f\,D_\pm s. 
\end{align}

Since one verifies easily that (cf. \cite[(4.2)]{MZ08})
\begin{align}\label{3}
\widetilde D_{f,\pm} (fs)-P_fD_\pm(fs)=-P_f\left(c(df)s\right), 
\end{align}
one sees that $\widetilde D_{f,\pm}$ is a compact perturbation of $P_fD_\pm$. Thus, one has
\begin{align}\label{4}
{\rm ind}\left(\widetilde D_{f,+}\right)={\rm ind}\left( P_fD_+\right).  
\end{align}

Now by (\ref{2}), if $fs\in\ker(\widetilde D_{f,+})$, then $s\in \ker (D_+)$.  Thus,   by the standard Lichnerowicz formula \cite{L63}, one has (cf. \cite[pp. 112]{GL83} and \cite[(3.6)]{L13})
\begin{align}\label{6}
 \frac{1}{2}\Delta\left(|s|^2\right)=\left|\nabla^{S_+(TM)} s\right|^2+\frac{k^{TM}}{4}|s|^2\geq \frac{k^{TM}}{4}|s|^2,
\end{align}
where $\Delta$ is the negative Laplace operator on $M$ and $\nabla^{S_+(TM)}$ is the canonical Hermitian connection on $S_+(TM)$ induced by $g^{TM}$. 

As has been observed in \cite{L13}, since the $G$-action on $M$ is cocompact and   $|s|$ is clearly $G$-invariant, there exists $x\in M$ such that 
\begin{align}\label{7a}
   |s(x)| ={\rm max}\{|s(y)|:{y\in M}\}. 
\end{align}
By the standard maximum principle, one has at $x$ that
\begin{align}\label{7}
 \Delta\left(|s|^2\right)\leq 0. 
\end{align}
Combining (\ref{7}) with (\ref{6}), one sees that if $k^{TM} > 0$ over $M$, one has 
\begin{align}\label{8}
 s(x)= 0, 
\end{align}
which implies that $s \equiv 0$ on $M$. 
Thus, one has $\ker(\widetilde D_{f,+})=\{0\}$, and, consequently, 
\begin{align}\label{9}
{\rm ind}\left(\widetilde D_{f,+}\right)\leq 0. 
\end{align}

On the other hand, for any $s,\ s'\in \Gamma(S(TM))$, one verifies that 
\begin{align}\label{10}
\left\langle fDs, fs'\right\rangle=\left\langle s, D\left(f^2s'\right) \right\rangle =\left\langle fs, D\left(fs'\right)+c(df)s'\right\rangle.
\end{align}

Let   $\widehat D_{f,\pm}: {\bf H}^1_f(M,S_\pm(TM) )^G \rightarrow {\bf H}^0_f(M,S_\mp(TM) )^G
$  be   defined by that for any $s\in \Gamma(S_\pm(TM) )^G$,
\begin{align}\label{11}
\widehat D_{f,\pm} (fs)=P_f\left(D_\pm(f s)+c(df)s\right) .
\end{align}
Clearly, $\widehat D_{f,+}$ is a compact perturbation of $P_fD_+$. Thus one has
\begin{align}\label{12}
{\rm ind}\left(\widehat D_{f,+}\right)={\rm ind}\left( P_fD_+\right).  
\end{align}

Now by (\ref{10}), one sees that the formal adjoint of $\widehat D_{f,+}$ is $\widetilde D_{f,-}$, while by proceeding as in (\ref{6})-(\ref{8}), one finds that
$\ker(\widetilde D_{f,-})=\{0\}$. Thus, one has
\begin{align}\label{13}
{\rm ind}\left(\widehat D_{f,+}\right)\geq 0.  
\end{align}

From (\ref{4}), (\ref{9}), (\ref{12}) and (\ref{13}), one gets
$
{\rm ind}\left( P_fD_+\right)= 0
$,
which completes the proof of Theorem \ref{t2}. 

$\ $

\noindent{\bf Acknowledgements.} The author is indebted to  Mathai Varghese for sharing his ideas in the joint work \cite{MZ08} and to Ziran Liu for helpful discussions.  {This work was
partially supported by MOEC and NNSFC.
}

\def\cprime{$'$} \def\cprime{$'$}
\providecommand{\bysame}{\leavevmode\hbox to3em{\hrulefill}\thinspace}
\providecommand{\MR}{\relax\ifhmode\unskip\space\fi MR }
\providecommand{\MRhref}[2]{%
  \href{http://www.ams.org/mathscinet-getitem?mr=#1}{#2}
}
\providecommand{\href}[2]{#2}

\end{document}